\documentclass[preprint,11pt]{elsarticle}

\usepackage[utf8]{inputenc}
\usepackage[english]{babel}
\usepackage{graphicx}
\usepackage{amsmath, accents}
\usepackage{amssymb}
\usepackage{amsfonts}
\usepackage{amsthm}
\usepackage{mathrsfs}
\usepackage{latexsym}
\usepackage{subfigure}
\usepackage{color}
\usepackage{xcolor}
\usepackage[noend]{algpseudocode}

\newcommand{\Z}{\mathbb{Z}}
\newcommand{\R}{\mathbb{R}}

\newtheorem{theorem}{Theorem}
\newtheorem{ex}[theorem]{Example}

\newtheorem{prop}[theorem]{Lemma}
\newtheorem{rem}[theorem]{Remark}
\newtheorem{cor}[theorem]{Corollary}

\journal{Discrete Applied Mathematics}

\begin{document}

\begin{frontmatter}

\title{Algorithms for linear time reconstruction by discrete tomography II}
\author[a] {Matthew Ceko}\ead{matthew.ceko@monash.edu}
\author[b]{Silvia M.C.~Pagani\corref{cor1}%
\fnref{fn1}}\ead{silvia.pagani@unicatt.it}
\author[c]{Rob Tijdeman}\ead{tijdeman@math.leidenuniv.nl}
\address[a] {School of Physics and Astronomy, Monash University, Melbourne, Australia}
\address[b]{Dipartimento di Matematica e Fisica ``N.~Tartaglia'', Università Cattolica del Sacro Cuore, via Musei 41, 25121 Brescia, Italy}
\address[c]{Mathematical Institute, Leiden University, 2300 RA Leiden, P.O.~Box 9512, The Netherlands}

\cortext[cor1]{Corresponding author}
\fntext[fn1]{Research supported by D1 Research line of Università Cattolica del Sacro Cuore.}

\begin{abstract}
The reconstruction of an unknown function $f$ from its line sums is the aim of discrete tomography. However, two main aspects prevent reconstruction from being an easy task. In general, many solutions are allowed due to the presence of the switching functions. Even when uniqueness conditions are available, results about the NP-hardness of reconstruction algorithms make their implementation inefficient when the values of $f$ are in certain sets. We show that this is not the case when $f$ takes values in a field or a unique factorization domain, such as $\R$ or $\Z$. We present a linear time reconstruction algorithm (in the number of directions and in the size of the grid), which outputs the original function values for all points outside of the switching domains. Freely chosen values are assigned to the other points, namely, those with ambiguities. Examples are provided.
\end{abstract}


\end{frontmatter}

\section{Introduction}
Tomography deals with the reconstruction of an object from the knowledge of its projections in a number of given directions. Radon \cite{radon} proved in 1917 that a differentiable function on $\mathbb{R}^2$ can be determined explicitly by means of integrals over the lines in $\mathbb{R}^2$. By approximating this for a large number of projections and using filtered back projection, so-called computerized tomography provides a quick way to compute a very good representation of the object. This method has a wide range of applications, from scans in hospitals to archaeology, astrophysics and industrial environments. See e.g.~\cite{her,nat}.

If the number of projection directions is small, discrete tomography may be advantageous compared to conventional back projection techniques. In this paper we consider a function $f$ on a finite grid $A$ of $\mathbb{Z}^2$ representing the object. Projections become line sums, i.e.~sums of the $f$-values at grid points on each line in finitely many given directions. Discrete tomography finds its origin in the fifties, mainly for only two directions, see e.g.~\cite{rys}. In 1978, Katz \cite{katz} gave a necessary and sufficient condition for the presence of a nontrivial function with vanishing line sums, known as a switching function or ghost. The theory started to blossom in the nineties when it became relevant in the study of crystals. In 1991 Fishburn, Lagarias, Reeds and Shepp \cite{fishepp} gave necessary and sufficient conditions for uniqueness of reconstruction of functions $f : A \to \{1,2, \ldots, N\}$ for some positive integer $N$.

An important distinction is whether the line sums are exact or may be inconsistent because of errors, termed noise, in the measurements. In case of noise the reconstruction can only be an approximation, see e.g.~\cite{bp,dart,pelt}. In what follows, we assume that the line sums are exact.

One of the main goals of discrete tomography is to ensure that the reconstructed function is equal to the function $f$ from which the line sums originate. However, in general the problem is ill-posed. Therefore one investigates which additional constraints can be imposed in order to achieve uniqueness. For instance, one may use some known information about the shape of the domain of $f$ such as convexity \cite{gg}, the values $f$ can attain (for the binary case see \cite{bdp1,hajdu}, for the integer case see \cite{dht}), or the size of the domain of $f$, \cite{bdp1,hatij}. In this paper we assume that the line sums come from some function $f$ and are therefore consistent. 

In 1999 Gardner, Gritzmann and Prangenberg \cite{ggp99} showed that the problem of reconstructing a function $f: A \to \mathbb{N}$ from its line sums in $d$ directions is solvable in polynomial time if $d=2$, but it is NP-complete if $d\geq 3$. The NP-completeness concerns both consistency and uniqueness, as well as reconstruction. Moreover, a year later they showed that the three mentioned problems are NP-complete for two and more directions when more than five types of atoms are involved in the crystal \cite{ggp00}.

We recall that the tomographic problem may be rephrased in terms of a linear system. If the function to be reconstructed has $\R$ as codomain, then it is known that polynomial-time algorithms exist to solve the linear system (such as the Gauss elimination, see \cite{atkinson}). The crux of the NP-results in \cite{ggp99,ggp00} is therefore the requirement that the range of $g$ is not closed under subtraction.


In 2001 Hajdu and Tijdeman \cite{hatij} gave an algebraic representation of the complete set of solutions over the integers. Their result also holds for solutions over the reals or any unique factorization domain. They gave a polynomial expression for the nontrivial switching function with domain of minimal size, the so-called primitive switching polynomial, and showed that every switching polynomial is a multiple of the primitive switching polynomial. This implies that every switching function is a linear combination of domain shifts of the corresponding primitive switching function. Their result implies that arbitrary function values can be given to a certain set of points and that thereafter the function values of the other points of $A$ are uniquely determined by the line sums. This was made explicit by Dulio and Pagani \cite{dupa} and serves as a building block in this paper.

In 2015 Dulio, Frosini and Pagani \cite{urrpdt,deda} showed that in the corners of $A$ the function values are uniquely determined and can be computed in linear time if the number of directions $d=2$. Later they proved conditional results for $d=3$ \cite{dgci2016,3dirext}. Recently, Pagani and Tijdeman \cite{pati} generalized the result for any number of directions. In particular the object function can be reconstructed in linear time if there are no switching functions. Moreover, they showed that in general the part of $A$ outside the convex hull of the union of all switching domains is uniquely determined and can be reconstructed in linear time. This result is another building block of our paper.


We prove that given the line sums of a function $f : A \to \R$ in the directions of a set $D$ we can compute a function $g: A \to \R$ with the same line sums. Using the theory of \cite{hatij} this implies that the complete set of such functions $g$ can be explicitly presented.

\begin{figure}
	\centering
	\setlength{\tabcolsep}{1pt}
	\begin{tabular}{cc}
		\includegraphics[width=0.49\columnwidth]{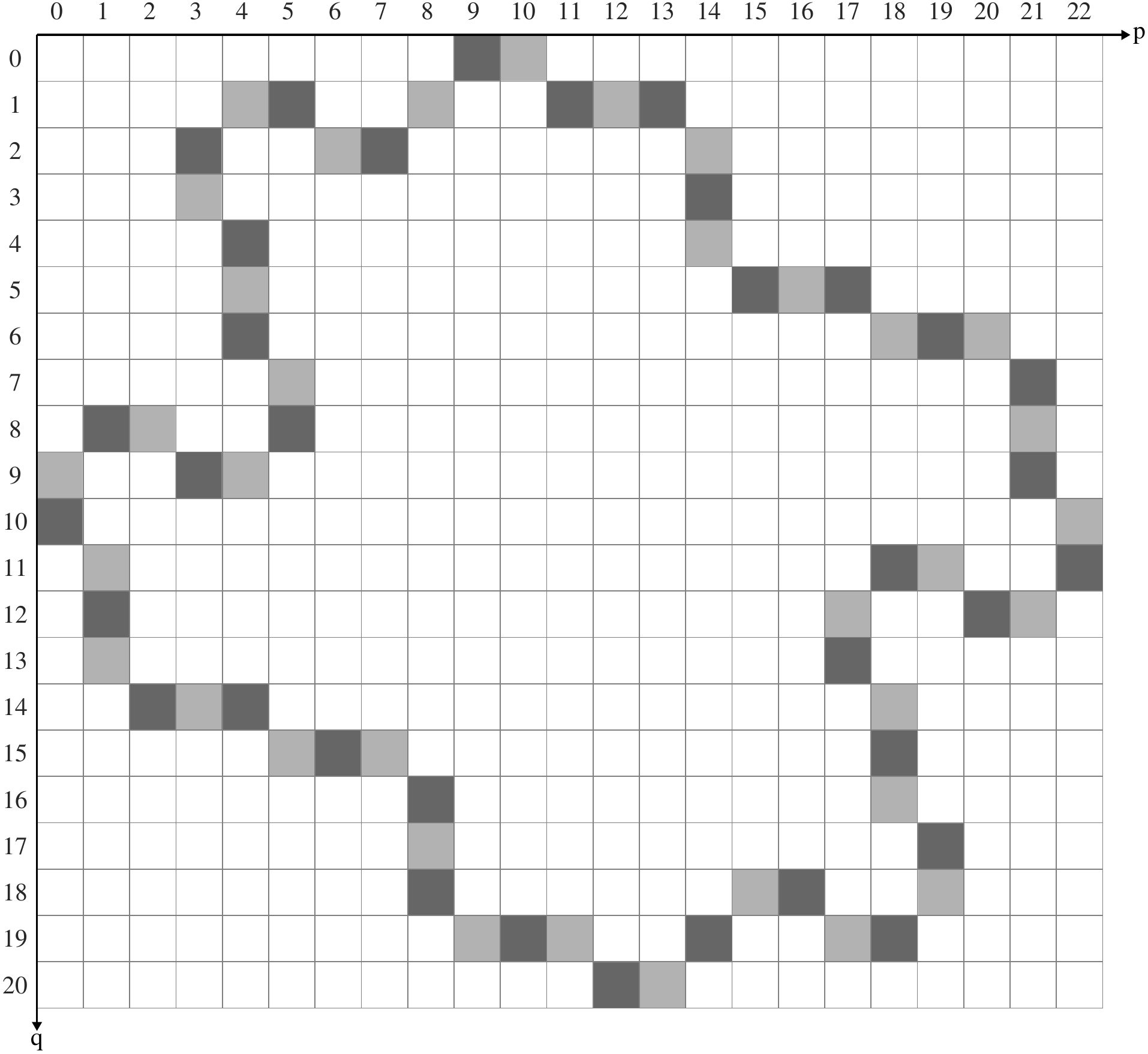}
	\end{tabular}
	\caption{Boundary ghost. The grey pixels form a switching domain. The pixels inside the ghost have $f$-values which are uniquely determined by the line sums in the directions of $D = \{(0,1), (1,0), (1,1), (-1,1), (-3,-1), (-1,-3), (5,-1), (7,5), (-3,7)\}$.}
	\label{fig_nearmin}
\end{figure}

Recently Ceko, Petersen, Svalbe and Tijdeman \cite{cpst} constructed switching components called boundary ghosts, where the switching domain has the form of an annulus around a relatively large interior, see e.g.~Figure \ref{fig_nearmin}. The values of $f$ for points which do not lie on this annulus can be uniquely determined by their line sums. This paper introduces a method which makes it possible to compute these values in linear time.


The present paper relies heavily on \cite{pati}, which was submitted before we started the research for the present paper. The above mentioned paper \cite{cpst} did us realize that it is important to be able to compute quickly the function values at the points in the interior of the boundary ghost domain. To our surprise we discovered that a twofold extension of the method of \cite{pati} worked, even for arbitrary ghosts. We explain this twofold application in the present paper. For our method it suffices that the range of $g$ is closed under subtraction. Further we provide a pseudo-code and a better justification of the linearity for the complexity than we did in \cite{pati}.

In Section \ref{sec:notation} we present notation and definitions, as well as information on switching functions. Section \ref{invalid} shows how values of $f$ in a corner region of $A$ can be obtained from the line sums. The case without switching components is treated in Section \ref{invalidcase}, that with switching components in Section \ref{sec:valid}. A general algorithm to compute $g$ can be found in Section \ref{effalg}. The justification of our linear time claim is given in Section \ref{sec:complexity}. Conclusions are in Section \ref{sec:conclusion}.


\section{Definitions and known results}\label{sec:notation}
We consider an $m\times n$ rectangular grid of points $$A=\{(p,q)\in\Z^2:0 \leq p< m,\,0 \leq q< n\}.$$ In our figures the $x$-axis is oriented from left to right and the $y$-axis from top to bottom. The origin is therefore the upper-left corner point of $A$. For each point $(p,q) \in \Z^2$ we consider the pixel $\{(x,y)\in\R^2\,:\,p \leq x<p+1, q \leq y <q+1\}$. In figures the coordinates of a pixel are the coordinates of the attached point.

Primitive directions are pairs $(a,b)$ of coprime integers. We agree to identify directions $(a,b)$ and $(-a,-b)$. Since we only consider primitive directions, we simply call them directions. The horizontal and the vertical direction are given by $(1,0)$ and $(0,1)$, respectively. We consider a finite set of directions $D= \{(a_h,b_h)~:~h=1,\ldots,d\}$. We say that $D$ is valid for $A$ if $ M :=\sum_{h=1}^d a_h <m$ and $N := \sum_{h=1}^d |b_h| <n$, and nonvalid otherwise.

A lattice line $L$ is a line containing at least two points in $\Z^2$. Let $f : A \to \R$. The line sum of $f$ along the lattice line $L(a,b,c): ay=bx+c$ with direction $(a,b)$ is defined as
$$\ell(a,b,c,f)=\sum_{aq=bp+c,~(p,q)\in A} f(p,q).$$

A function $F: A \to \mathbb{R}$ is called a \emph{switching function} or \emph{ghost} of $(A,D)$ if all the line sums of $F$ in all the directions of $D$ are zero. Observe that then $f$ and $f+F$ have the same line sums in the directions of $D$. The support of a switching function is called a \emph{switching domain}.

We say that something can be computed in linear time if the number of basic operations needed to compute it is $\mathcal{O}(dmn)$. Here a basic operation is an addition, subtraction, multiplication, division, decision about which of two quantities is larger or an assignment.

\subsection{The location of switching domains} \label{sec_location}

M.~Katz \cite{katz} proved that $f: A \to \R$ is uniquely determined by the line sums in the directions of $D$ if and only if $(A,D)$ is nonvalid. Fishburn et al.~\cite{fishepp} showed that $(p,q) \in A$ has a unique $f$-value if and only if $(p,q)$ is not located in a switching domain. Hajdu and Tijdeman \cite{hatij} associated to the function $f: A \to \mathbb{R}$ the polynomial $f^*(x,y) = \sum_{i=0}^{m-1} \sum_{j=0}^{n-1} f(i,j) x^iy^j$. In this way every switching function corresponds with a switching polynomial. They defined
\begin{displaymath}
g^*_{(a,b)}(x,y) = \left\{
\begin{array}{ll}
x^ay^b-1 & \text{ if } a>0, b>0,\\
x^a-y^{-b} & \text{ if } a>0, b<0,\\
x-1 & \text{ if } a=1, b=0,\\
y-1 & \text{ if } a=0, b=1,
\end{array} \right.
\end{displaymath}
and $$G^*_{i,j}(x,y)=x^iy^j \prod_{h=1}^d g^*_{(a_h,b_h)}(x,y)$$ for $0 \leq i < m- M, 0 \leq j < n - N$. They showed that $G^*_{0,0}$ is a switching polynomial of minimal degree. We call the corresponding function a primitive switching function. Furthermore they proved the following result.

\begin{theorem} [Hajdu, Tijdeman \cite{hatij}, Theorem 1] \label{hati}
Suppose $D$ is valid for $A$. Put $M= \sum_{h=1}^d a_h, N= \sum_{h=1}^d |b_h|.$ Then for every switching function $g: A \to \R$ its switching polynomial $g^*$ can be uniquely written as
\begin{equation} \label{swfu}
g^* = \sum_{i=0}^{m-1-M} ~\sum_{j=0}^{n-1-N} c_{i,j} G^*_{i,j}
\end{equation}
with $c_{i,j} \in \R$ for all $i,j$. Conversely, every function $g$ of which the switching polynomial is of the form \eqref{swfu} is a switching function.
\end{theorem}

\begin{figure}[htbp]
\centering
\includegraphics[scale=0.47]{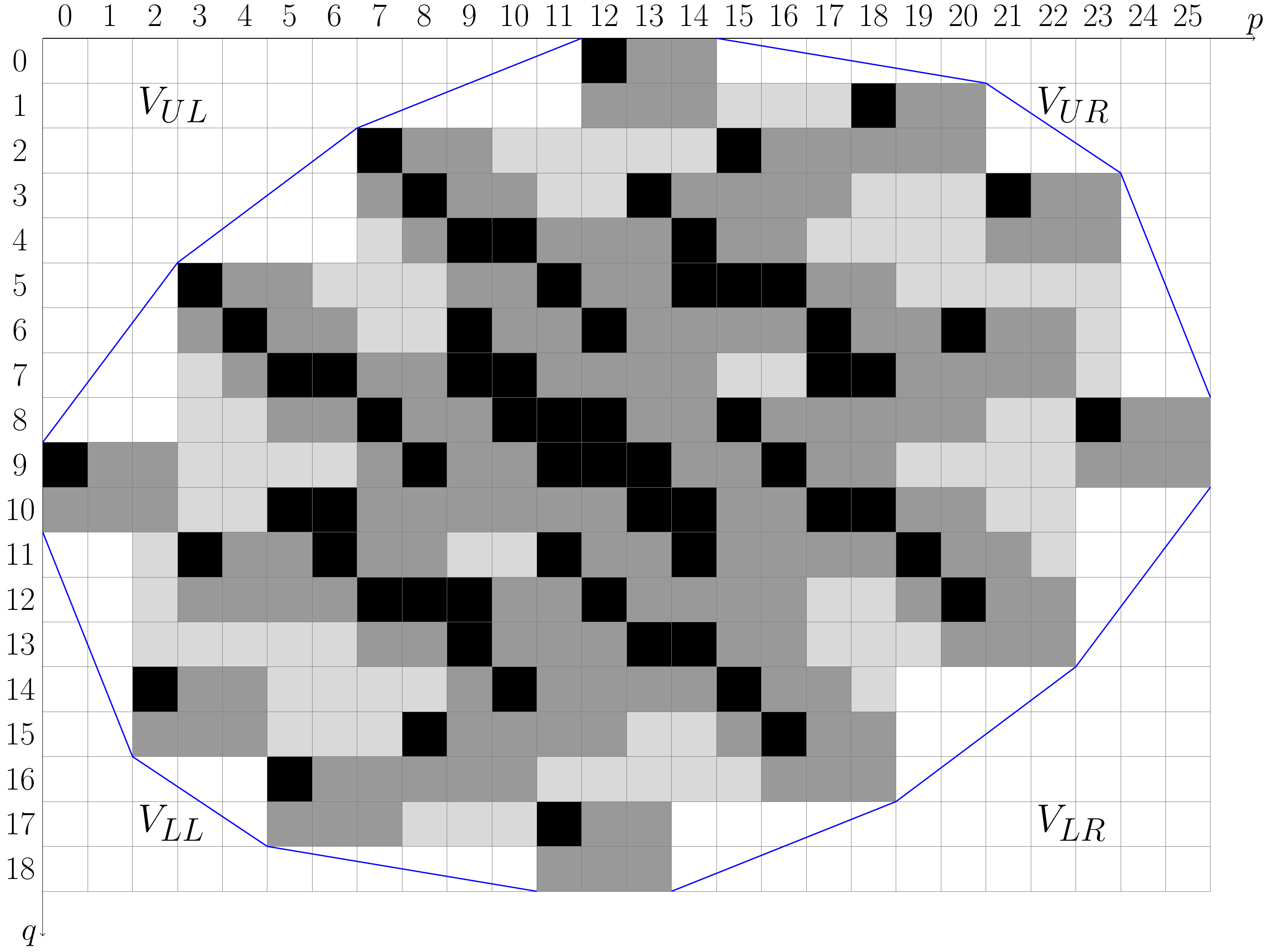}
\caption{\small The situation for the 26 by 19 grid $A$ and the set of directions $D=\{(5,-2),(4,-3),(3,-4), (6,1), (3,2), (2,5)\}$. The dark grey and black pixels indicate the union of the switching domains. The black pixels represent the switching domain related to $G^*_{0,0}$. In the dark grey and black pixels the function $f$ is not uniquely determined by its line sums in the directions of $D$. The $f$-values of the complement, the white and light grey pixels, are uniquely determined by these line sums. Since $M=23, N=17$, there are six pixels where the choice is free, e.g.~the pixels $(0,9)$, $(0,10)$, $(1,9)$, $(1,10)$, $(2,9)$, $(2,10)$. Any other 3 by 2 block of dark grey and black pixels can be chosen instead. If the choice is made all the values of the unique solution satisfying the made choices are determined by the line sums in the directions of $D$. The white pixels form four corner regions. The broken line indicates the convex hull of the union of the switching components.}
\label{fig:union}
\end{figure}

This result is also valid if $\R$ is replaced by $\Z$ or any other field or unique factorization domain. A corollary of the theorem relevant for this paper is that the lexicographically lowest degree term of $G_{i,j}^*$ is given by $x^iy^{j+N_n}$ where $N_n = \sum_{b_h<0} -b_h$. Thus we have free choice for the values of $c_{i,j}$ for $0 \leq i < m-M, N_n \leq j < N_n+n-N$ and by this choice the function $g^*$ is uniquely determined. An illustration of Theorem \ref{hati} is given in Figure \ref{fig:union}.


\section{Uniqueness in the corner regions} \label{invalid}
\begin{figure}[htbp]
\centering
\includegraphics[scale=0.7]{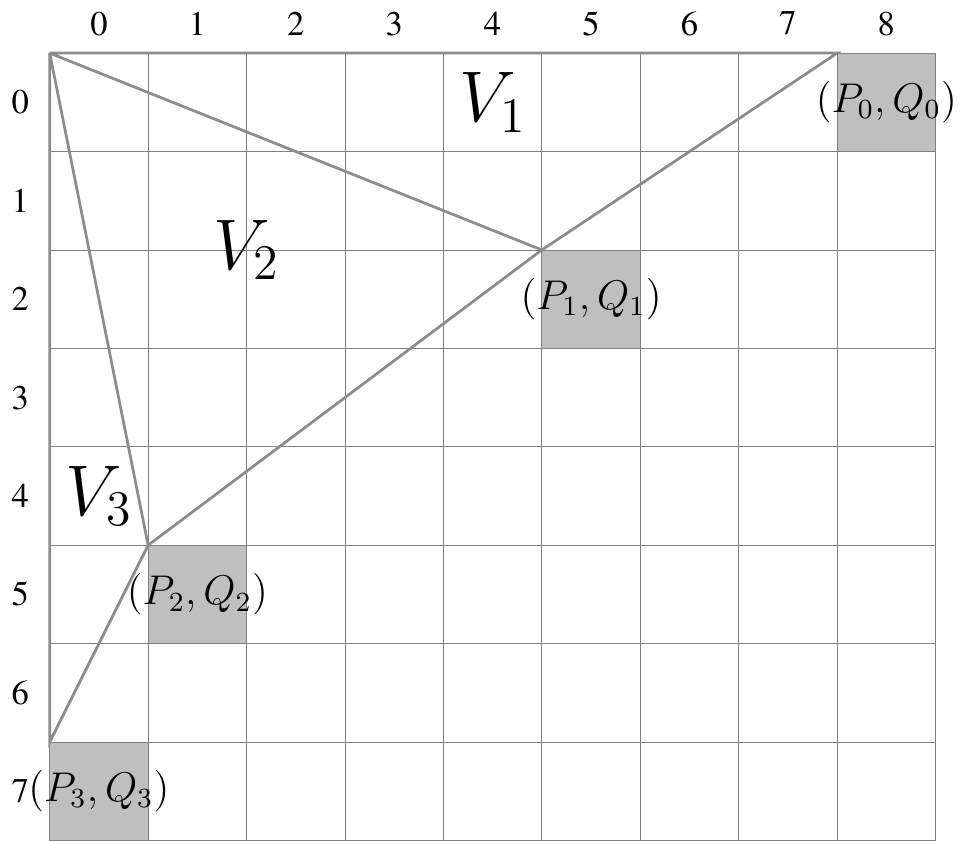}
\caption{\small The triangles $V_1,V_2,V_3$ for the set $D=\{(3,-2),(4,-3),(1,-2)\}$. The border points are $(P_0,Q_0)=(8,0),(P_1,Q_1)=(5,2),(P_2,Q_2)=(1,5),(P_3,Q_3)=(0,7)$. For every $h$ the line through $\left(P_{h-1},Q_{h-1}\right)$ and $\left(P_h,Q_h\right)$ is a side of triangle $V_h$, and intersects each other triangle $V_{\tilde{h}}$, since the slopes increase with increasing $\tilde{h}$ by the ordering in \eqref{ratios}.}
\label{fig:triangles1}
\end{figure}

Let again $A=\{(p,q)\in\Z^2:0 \leq p< m,\,0 \leq q< n\}$. Let $D$ be a set of directions $(a_1, -b_1), \ldots, (a_k,-b_k)$ with $k \geq 2$ where $a_1, \ldots, a_k, b_1, \ldots b_k$ are positive integers ordered such that
\begin{equation}\label{ratios}
\frac{b_1}{a_1}<\frac{b_2}{a_2} < \ldots < \frac {b_k}{a_k}.
\end{equation}
Note that by primitivity all the ratios are distinct. We call the points
$$\left( \sum_{h=1}^k a_h,\, 0\right),\, \left(\sum_{h=2}^ka_h,\, b_1\right),\, \left(\sum_{h=3}^ka_h,\, \sum_{h=1}^2b_h\right),\, \ldots,\, \left(0,\, \sum_{h=1}^k b_h\right)$$
the border points of the upper left region $(P_0,Q_0), (P_1,Q_1), \ldots, (P_k,Q_k)$, respectively. We denote the convex hull of the three points $(0,0)$, $(P_{h-1},Q_{h-1})$, $(P_h,Q_h)$ by $V_h$ for $h=1,2, \dots, k$ (see Figure \ref{fig:triangles1}). Let
$$V_{UL} = \bigcup_{h=1}^k V_h$$
be the upper left corner region. The other corner regions $V_{UR}$, $V_{LL}$, $V_{LR}$ may be defined similarly (see Figure \ref{fig:union}).

For a point $(p,q) \in A$ we define its weight $w(p,q)$ by
$$w(p,q) = \min_{h=1,2,\ldots,k} ~~\frac {b_hp+a_hq}{b_h P_h+ a_h Q_h}.$$
The weight function in $V_{UL}$ equals the quotient of the distance of the point $(p,q)$ to the origin $(0,0)$ and the distance from the origin to the intersection $(p',q')$ of the line through $(0,0)$ and $(p,q)$ and the boundary of the convex hull. This weight has the property that every point $(p,q)$ in $V_{UL}$ has maximal weight among the integer points on the line $\ell$ through $(p,q)$ parallel to the line segment of the boundary of the convex hull through $(p',q')$.

\begin{figure}[htbp]
\centering
\includegraphics[scale=0.7]{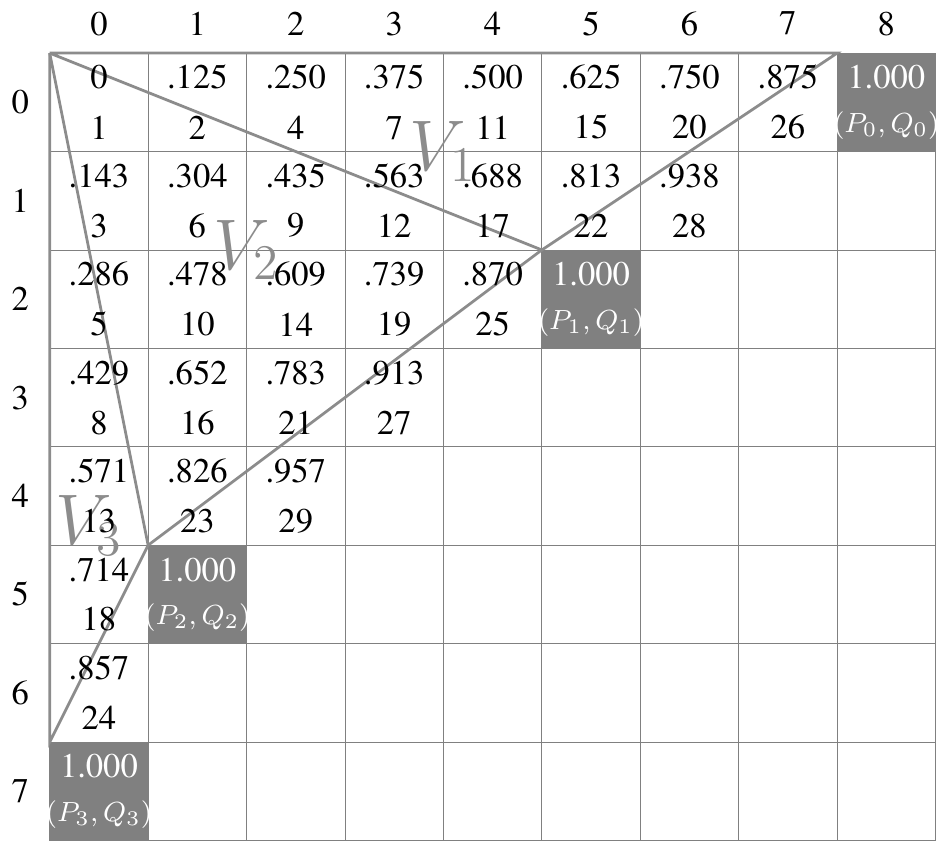}
\caption{\small The weights (upper number inside each pixel) of the points for directions $(3,-2), (4,-3), (1,-2)$. The border points are $(8,0),(5,2),(1,5), (0,7)$. All the points with weight less than $1$ are in the corner region and have uniquely determined $f$-values (Theorem \ref{aROU}). The lower numbers enumerate them with increasing weights. The (dark grey) border pixels are part of the switching domain and their $f$-values are therefore not uniquely determined. They have weight $1$. All entirely white pixels have weight $>1$.}
\label{fig:weight}
\end{figure}

The following lemma implies that if $(p,q) \in V_{UL}$, then the minimum in the definition of $w(p,q)$ is reached for $h$ such that $(p,q) \in V_h$ (see Figure \ref{fig:weight}).

\begin{prop}[\cite{pati}, Lemma 2] \label{weight}
For $(p,q) \in A$ the weight $w(p,q)$ is reached for $h$ such that
\begin{displaymath} 
\frac{Q_{h-1}}{P_{h-1}} \leq \frac qp \leq \frac{Q_h}{P_h}.
\end{displaymath}
and only for such $h$. The weight $1$ is reached at the border points and not at other points of $A$.
\end{prop}

The next result states that the corner region $V_{UL}$ except for the border points has unique $f$-values which can be computed in linear time.

\begin{theorem}[\cite{pati}, Theorem 4 and Corollary 6] \label{aROU}
Let $A=\{(p,q)\in\Z^2:0 \leq p< m,\,0 \leq q< n\}$. Let $D$ be a set of directions $(a_1, -b_1), \ldots, (a_k,-b_k)$ where $a_1, \ldots, a_k, b_1, \ldots b_k$ are positive integers ordered as in \eqref{ratios}. Let the line sums of $f : A \to \R$ in the directions of $D$ be given. Then all the points $(p,q)$ in $V_{UL}$ except for the border points have uniquely determined $f$-values.
\end{theorem}

\begin{cor}[\cite{pati}]
The values of the points as in the previous theorem can be computed according to increasing weights and, if $(p,q) \in V_h$, by subtracting the sum of the $f$-values of the other points of $A$ on the line through $(p,q)$ in the direction of $(a_h,b_h)$ from its line sum.
\end{cor}

\begin{cor} \label{cor_column}
Under the above conditions the $f$-values of the points $(0,0)$, $(0,1)$, \ldots, $\left(0, -1+\sum_{h=1}^k b_h\right)$ can all be computed in linear time.
\end{cor}


\section{The nonvalid case} \label{invalidcase}

Suppose we are in the nonvalid case, then $M \geq m$ or $N \geq n$. Without loss of generality assume that $N \geq n$. Then we apply Theorem \ref{aROU} both to the upper corner region $V_{UL}$ and to the lower corner region $V_{LL}$.

Let $A$ be as above. Let $$D =\{(a_1, -b_1), \ldots, (a_k,-b_k), (a_{k+1}, b_{k+1}), \ldots, (a_d,b_d), (0,1)^*, (1,0)^*\}$$ where $a_1, \ldots, a_d, b_1, \ldots b_d$ are positive integers ordered such that
\begin{displaymath}
\frac{b_1}{a_1}<\frac{b_2}{a_2} < \ldots < \frac {b_k}{a_k},\qquad\frac{b_{k+1}}{a_{k+1}}> \frac{b_{k+2}}{a_{k+2}}>\ldots > \frac {b_d}{a_d}
\end{displaymath}
and the asterisk indicates that $(0,1)$ and/or $(1,0)$ may occur in $D$. Thus we assume that $n \leq \sum_{h=1}^d b_h$ or ($n = 1+\sum_{h=1}^d b_h$ and $(0,1) \in D)$.

By Corollary \ref{cor_column} applied to $V_{UL}$, the $f$-values of the points $(0,0)$, $(0,1)$, \dots, $\left(0, -1+\sum_{h=1}^k b_k\right)$ can be computed. In a similar way we can apply the corollary to $V_{LL}$ and the directions $(a_{k+1}, b_{k+1}), \ldots, (a_d,b_d)$ to conclude that the $f$-values of the points $(0, n-1), (0,n-2), \dots, \left(0, n- \sum_{h=k+1}^d b_h\right)$ can be computed. It follows that the $f$-values of the points $(0,0), (0,1), \dots, (0,n-1)$ can all be computed except when $n = 1+\sum_{h=1}^d b_h$ and $(0,1) \in D$. In the latter case $(p,q)=\left(0, \sum_{h=1}^k b_h\right)$ is the only point in the column $p=0$ with unknown $f$-value. However, this value can be found by subtracting from the line sum of the column $p=0$ the $f$-values of the other points in that column. In this way we have made our problem of computing the $f$-values one column smaller. We can repeat the procedure in order to find the $f$ values of the next column. Continuing the process we arrive at the following conclusion.

\begin{theorem}[\cite{pati}]\label{thm:invalid}
Let $A=\{(p,q)\in\Z^2:0 \leq p< m,\,0 \leq q< n\}$. Let $D$ be a set of directions such that $A$ is nonvalid for $D$. Let the line sums of $f : A \to \R$ be given. Then the $f$-values of all points of $A$ can be computed in linear time.
\end{theorem}

In \cite{pati} algorithms are given for computing the $f$-values. These algorithms are more efficient than the procedure described above. The algorithm in Section \ref{effalg} is as efficient as these algorithms.


\section{The valid case}\label{sec:valid}

Let $A=\{(p,q)\in\Z^2:0 \leq p< m,\,0 \leq q< n\}$ and $$D =\{(a_1, -b_1), \ldots, (a_k,-b_k), (a_{k+1}, b_{k+1}), \ldots, (a_d,b_d), (0,1)^*, (1,0)^*\}$$ where $a_1, \ldots, a_d, b_1, \ldots b_d$ are positive integers ordered such that
\begin{displaymath}
\frac{b_1}{a_1}<\frac{b_2}{a_2} < \ldots < \frac {b_k}{a_k},\qquad\frac{b_{k+1}}{a_{k+1}}> \frac{b_{k+2}}{a_{k+2}}>\ldots > \frac {b_d}{a_d}
\end{displaymath}
and the asterisk indicates that $(0,1)$ and/or $(1,0)$ may occur in $D$. As observed in the previous section, by applying Corollary \ref{cor_column} to $V_{UL}$ the $g$-values of the points $(0,0), (0,1), \dots, \left(0, -1+\sum_{h=1}^k b_h\right)$ can be computed. In a similar way we can apply the corollary to $V_{LL}$ and the directions $(a_{k+1}, b_{k+1}), \ldots, (a_d,b_d)$ to conclude that the $g$-values of the points $(0, n-1), (0,n-2), \dots, \left(0, n- \sum_{h=k+1}^d b_h\right)$ can be computed. In Section \ref{sec_location} it was observed that the $g$-values of the points $(0, Q_k), (0, Q_k+1), \dots, (0, Q_k+n-N-1)$ can be freely chosen where $g : A \to \R$ is a function satisfying the line sums. Combining these results we see that all the $g$-values of the points $(0,0), (0,1), \dots, (0,n-1)$ can be computed or freely chosen, except for the case that $(0,1) \in D$ and $n=1+ \sum_{h-1}^d b_d$. In the latter case only the $g$-value of $\left(0, \sum_{h=1}^k b_h\right)$ is not determined, but this can be computed by subtracting the $g$-values of the other points in the leftmost column from the sum of that column. After this all $g$-values of the points in the leftmost column are fixed. In Section \ref{sec_location} it was further observed that the $g$-values of the points $(p, Q_k), (p, Q_k+1), \dots, (p, Q_k+n-N-1)$ for $p=1,2, \dots, m-M-1$ can be freely chosen. Therefore we can repeat the above procedure successively for columns $p=1,2, \dots, m-M-1$. Then $M$ columns remain, for which the line sums in the directions of $D$ are known. It is obvious that the computed $g$-values of the points which do not belong to a switching domain have the original $f$-value. We are left with a nonvalid case and we can apply an algorithm for that case to compute the remaining $g$-values.

We have shown that the following theorem holds.

\begin{theorem}\label{lineartime}
Let $A=\{(p,q)\in\Z^2:0 \leq p< m,\,0 \leq q< n\}$ and $D$ a set of primitive directions. Let $f : A \to \R$ be an unknown function. Suppose all the line sums in the directions of $D$ are known. Then we can compute a function $g : A \to \R$ satisfying the line sums in linear time. The points which do not belong to any switching domain get their original $f$-value.
\end{theorem}

\begin{ex}{\rm
Consider the situation in Figure \ref{fig:union}. We have $m=26, M=23, n=19, N=17, Q_k=9.$ We can freely choose the $g$-values of the points $(0,9)$ and $(0,10)$ and compute the $g$-values of the other points $(0,q)$. Next we do so for the columns $p=1$ and $p=2$. We are left with a 23 by 19 rectangular grid. Since $m=M=23$, this is a nonvalid case and we know that the remaining $g$-values can be computed in linear time. The found $g$-values of the white and light grey pixels are equal to the original $f$-values.}
\end{ex}

\begin{rem}{\rm
In this paper we assume that the line sums are correct and that there is no noise. It is easy to check whether this is true afterwards by checking the line sums which have not been used for computing the $g$-values. In case the line sums are inconsistent, and it is better to use a method which treats the unused line sums in a similar way as the used line sums to obtain a good approximation of the original function.}
\end{rem}

\begin{rem}{\rm
Theorem \ref{hati} states that if $g : A \to \R$ has the same line sums as $f$, then the associated polynomial $g^*$ is of the form
$$f^* + \sum_{i=0}^{m-1-M} ~\sum_{j=0}^{n-1-N} c_{i,j} G^*_{i,j}$$
and each such function has the same line sums as $f$. It is possible to compute the coefficients $c_{i,j}$ as follows. The point $(0,Q_k)$ occurs only in the domain of $G_{0,0}$ and therefore $c_{0,0}$ can be found from the found value for $(0,Q_k)$. The point $(0,Q_k+1)$ occurs in $G_{0,1}$ and maybe in $G_{0,0}$. Since $c_{0,0}$ is already known, $c_{0,1}$ can be computed. Considering the points with a free choice in the lexicographic order, each time a point occurs in only one new primitive switching domain and hence the corresponding coefficient can be computed.}
\end{rem}

\section{An efficient algorithm} \label{effalg}
In this section we present an algorithm to find a function $g : A \to \R$ which satisfies the given line sums of an unknown function $f : A \to \R$. This algorithm is based on ideas and results in \cite{pati}. Its complexity is studied in the next section.

In this algorithm it is not necessary to compute all weights as in Figure \ref{fig:weight}. Observe that after the $g$-values in $V_{UL}$ and $V_{LL}$ have been computed, the $g$-value of the next point on each row will be computed. For this the same direction will be used as used for the integer point immediately left of it, since everything shifts one place to the right. For the same reason the order in which the new points will be handled will be the same as for the points immediately left of them. This process will be continued. Thus, where in the example of Figure \ref{fig:weight} initially on the row $q=1$ first the direction $(1,-2)$, next the direction $(4,-3)$, and finally the direction $(3,-2)$ was used to compute the $g$-value, more to the right only the direction $(3,-2)$ would be used. Suppose there would have been seven more columns $p=-1, -2, \ldots, -7$ with $g$-values equal to 0 in $A$ (cf.~Algorithm 1B of \cite{pati}). Then for $p\geq0$ we would have needed only the rightmost direction on each row and the direction needed to compute the $g$-values of points in the original $A$ would only depend on their row. Therefore it suffices to follow the order in which the rightmost non-border points of $V_{UL}$ and $V_{LL}$ are treated and for each point $(p,q)$ to use the line sum in the direction $(a_h,b_h)$ with $h$ such that the corresponding rightmost point is in $V_h$. Observe that this $h$ is such that $Q_h \leq q < Q_{h+1}$.

\begin{figure}[htbp]
\centering
\includegraphics[scale=0.7]{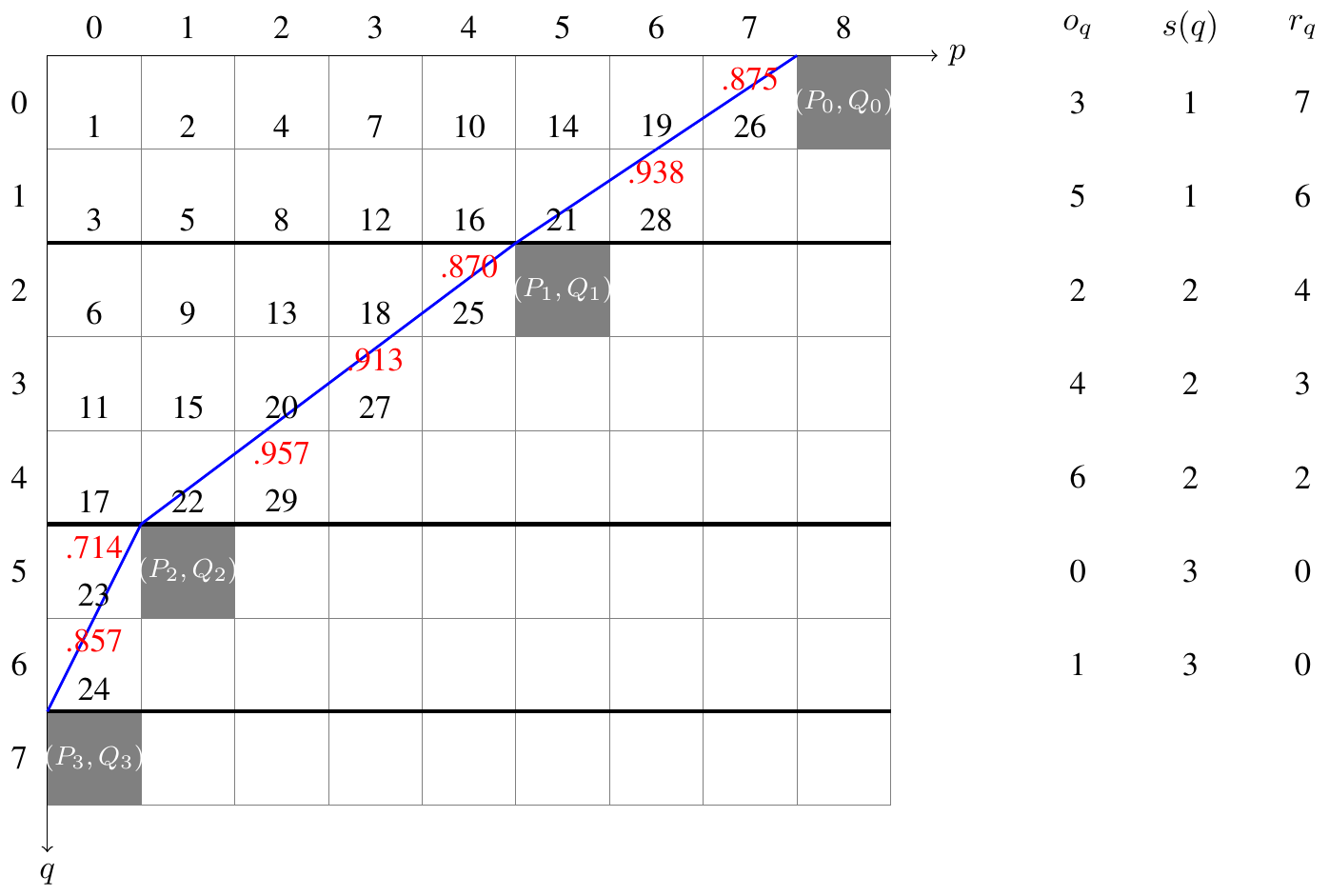}
\caption{\small The weights (upper numbers) of the rightmost points in $V_{UL}$ of each row and the ordering of points in $V_{UL}$ (lower numbers) of Example \ref{ex_neworder}. On the right, $o_q$ reports the order of the rightmost points in $V_{UL}$ after increasing weights, $s(q)$ tells which direction has been used in the corresponding row and $r_q$ indicates, for each row, the $p$-coordinate of the rightmost point in $V_{UL}$.}
\label{fig:highestweight}
\end{figure}

\begin{ex}\label{ex_neworder}{\rm
Let as in Figure \ref{fig:weight} the directions be $(3,-2), (4,-3), (1,-2)$. The border points are $(8,0),(5,2),(1,5), (0,7)$. Let $m >8, n=7.$ Then the rightmost points with weight $<1$ are $$(7,0),~ (6,1), ~(4,2), ~(3,3), ~(2,4),~ (0,5),~ (0,6).$$ If we order them according to increasing weights, then we get
$$(0,5), ~(0,6),~(4,2),~(7,0),~(3,3),~(6,1),~(2,4).$$
If we use the invisible seven columns on the left with $g$-values 0, then the enumeration of $V_{UL}$ is given by the lower numbers in Figure \ref{fig:highestweight}. Observe that it differs from the enumeration in Figure \ref{fig:weight}. For computing the $g$-values in rows 0 and 1 direction $(3,-2)$ is used, in rows 2 to 4 direction $(4,-3)$ and in rows 5 and 6 direction $(1,-2)$. 
}
\end{ex}

The above example illustrates the first seven steps of the algorithm. The above procedure is done for both $V_{UL}$ and for $V_{LL}$. In between $g$-values 0 are substituted (or any other values) at the places where the switching domains offer free choice. Now the $g$-values in the first column are known and we can proceed with the next column and so on until we have treated $m-M$ columns. An $M$ by $n$ grid remains to be handled, but this is a nonvalid case. This can be treated in a similar way, but starting from upper corner regions $V_{UL}$ and $V_{UR}$ and then going downwards. In the algorithm we have $g(p,q)=f(p,q)$ for all the points $(p,q)$ which are not in a switching domain of $(A,D)$.

\begin{figure}[htbp]
\centering
\includegraphics[scale=0.6]{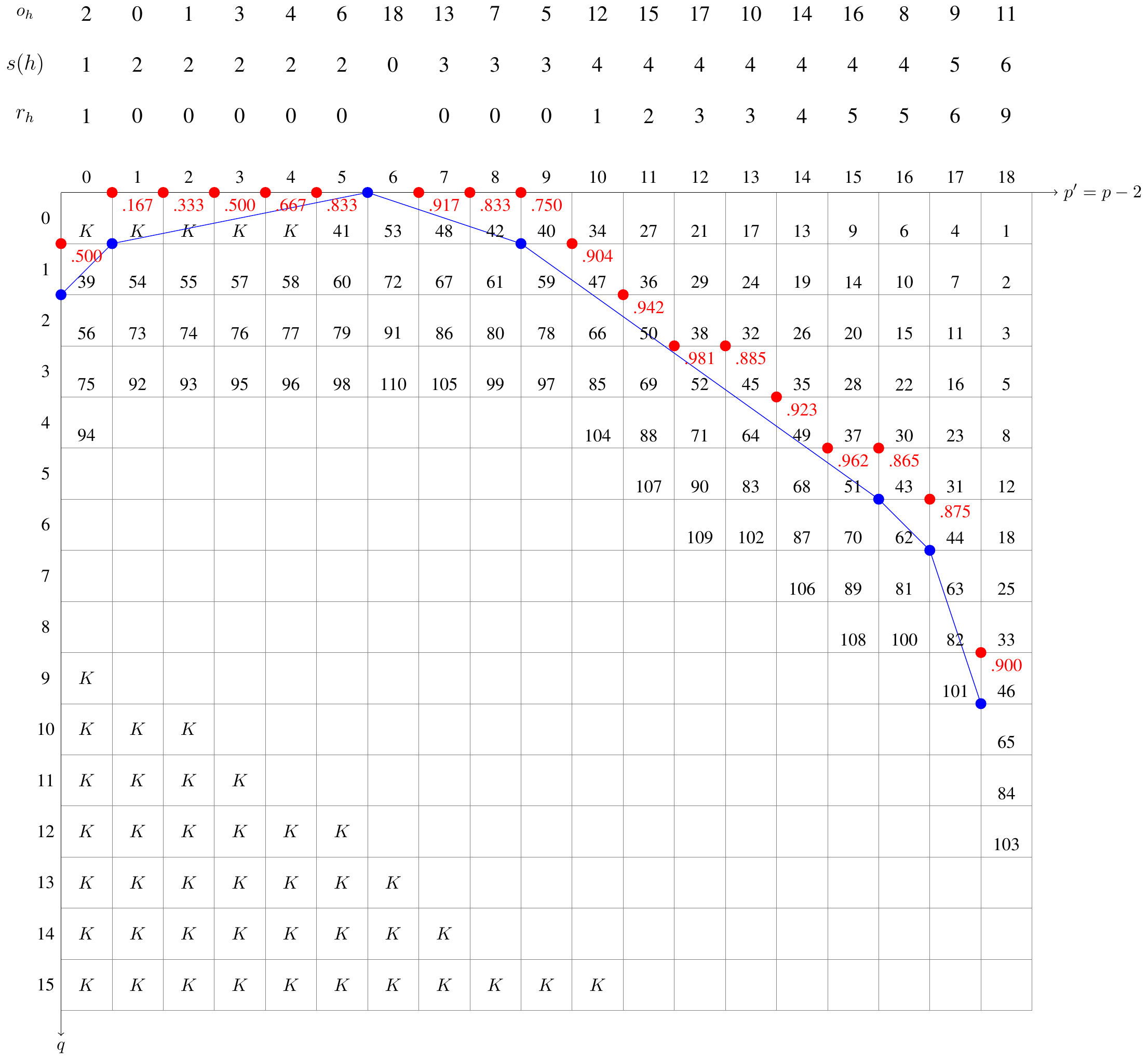}
\caption{\small This figure illustrates Steps 8-12 of Example \ref{ex_steps}. We have omitted the first $m-M=2$ columns so that the column numbers indicate $p-2$. The $K$'s indicate the pixels of which the $g$-values have already been calculated in Step 7. In this stage the roles of rows and columns are interchanged. Step 8 serves to adjust the order of the directions. In Step 9 the border points in $V_{UL}$ are reordered, those in $V_{UR}$ are those of $V_{LL}$ mirrored. Again the broken line connects them. This time the points which determine the direction to be used are the integer points immediately {\it above} the border points. Their $Q$-values, their weights (the above number inside the pixel) and the sequence $s$ are computed in Step 10. They are ordered in Step 11 and the order is indicated in row $o_h$. Finally in Step 12 the remaining $g$-values are computed where by using $u$ the original $p$-values are used instead of $p-2$. The order of the way the $g$-values are found is indicated by the lower numbers inside the pixels. This has to be completed downwards to find all $g$-values.}
\label{fig:algnonvalid}
\end{figure}

The algorithm has 12 steps and is illustrated in Example \ref{ex_steps} following it. For Steps 1-7 see Figure \ref{fig:algvalid}, for Steps 8-12 see Figure \ref{fig:algnonvalid}.

In Step 1 the directions are ordered in such a way that the corner regions become concave. Throughout the algorithm we write $b_h$ instead of $-b_h$ ($h=1,\dots,k$) so that $b_h$ is always nonnegative. In Step 2 the border points of $V_{UL}$ and $V_{LL}$ are found. In Step 3 a function $\delta$ is introduced indicating whether the $g$-value of a point has been computed. Step 4 provides a shortcut for nonvalid cases. If $n \leq N$, then this shortcut can be used after rows and columns have been interchanged. Step 5 serves to find the grid points left of the border line, measured by the sequence $r$. The weights of these points are computed as well, together with the sequence $s$, which indicates the direction of the line used to compute the $g$-value. In Step 6 the points found in the previous step are ordered after increasing weight by the sequence $o$. If weights are equal, the order is irrelevant. In Step 7 the $g$-values of the first $m-M$ columns (and of some more points of $A$) are computed.

Steps 8-12 are essentially equal to Steps 1-7, but with the columns $p=0,1, \dots, m-M-1$ omitted, as they have already been treated, and the roles of the rows and columns interchanged. In Step 8 the directions are reordered as now $V_{UL}$ is mirrored and $V_{UR}$ takes over the role of $V_{LL}$. A similar reordering of border points takes place in Step 9. The weights and the corresponding directions are found in Step 10. The order of the points found in the previous step are fixed in Step 11. Finally the remaining $g$-values are computed in Step 12 where the $u$ is introduced to make the necessary shift because of the omitted $m-M$ columns.

In the algorithm, $\lambda(\mathbf{d},p,q)$ denotes the line sum containing the point $(p,q)$ in the direction $\mathbf{d}$, and $\left\lceil r \right\rceil$ denotes the ceiling of $r$.


\vspace{.3cm}
\noindent \textbf{Algorithm.}
\begin{algorithmic}
    \Require{A set $A = \{(p,q) : 0 \leq p <m, 0 \leq q<n \}$ with positive integers $m,n$, a finite set of (primitive) directions $D$ and all the line sums in the directions of $D$ of a function $f : A \to \R.$}
    \Ensure{Function $g : A \to \R$ which satisfies the line sums.}
\end{algorithmic}
{\bf Step 1:} Initial values.
\begin{algorithmic}
    \renewcommand\algorithmicdo{}
    \ForAll{$\mathbf{d}=(a,-b) \in D$ (with $a>0,b>0$)}order the directions such that
        \begin{equation*}
            \frac{b_1}{a_1}<\frac{b_2}{a_2} < \dots < \frac {b_k}{a_k}.
        \end{equation*}
    \EndFor
    \ForAll{$\mathbf{d}=(a,b) \in D$ (with $a>0,b>0$)}order the directions such that
        \begin{equation*}
            \frac{b_{k+1}}{a_{k+1}}>\frac{b_{k+2}}{a_{k+2}} >\dots > \frac {b_d}{a_d}.
        \end{equation*}
    \EndFor
    \renewcommand\algorithmicdo{\textbf{do}}
    \State $M\gets\sum_{h=1}^d a_h$
    \State $N\gets\sum_{h=1}^d b_h$
    \If {$(1,0) \in D$} $M \gets M+1$ \EndIf
    \If {$(0,1) \in D$}
        \State $N \gets N+1$
        \State $\mathbf{d}_0 \gets (0,1)$
    \EndIf
\end{algorithmic}
{\bf Step 2:} Border points.
\begin{algorithmic}
    \State $(P_0,Q_0) \gets \left(\sum_{h=1}^k a_h,0\right)$
    \For{$h \gets 1$ {\bf to} $k$}
        \State $(P_h,Q_h) \gets (P_{h-1}-a_h, Q_{h-1}+b_h)$
    \EndFor
    \State $(P^*_k,Q^*_k) \gets \left(0, n-1-\sum_{j=k+1}^d b_j\right)$
    \State $(P_{k+1},Q_{k+1}) \gets (P^*_k +a_{k+1}, Q^*_k+b_{k+1})$
    \For{$h \gets k+2$ {\bf to} $d$}
        \State $(P_h,Q_h) \gets (P_{h-1}+a_h,Q_{h-1}+b_h)$
    \EndFor
\end{algorithmic}
{\bf Step 3:} Fixing switching functions.
\begin{algorithmic}
    \For{$p \gets 0$ {\bf to} $m-1$}
        \For{$q \gets 0$ {\bf to} $n-1$}
            \State $\delta(p,q) \gets 0$
        \EndFor
    \EndFor

    \For{$p \gets 0$ {\bf to} $m-M-1$}
        \For{$q \gets Q_k$ {\bf to} $Q_k+n-N-1$}
            \State $g(p,q) \gets 0$
            \State $\delta(p,q) \gets 1$
        \EndFor
    \EndFor
\end{algorithmic}
{\bf Step 4:} Nonvalid case.
\begin{algorithmic}
    \If{$m \leq M$} goto Step 8\EndIf
\end{algorithmic}
{\bf Step 5:} Choosing starting points $r_h$, weights $w(r_h,h)$ and directions $\mathbf{d}_{s(h)}$.
\begin{algorithmic}
    \For{$H \gets 1$ {\bf to} $k$}
        \For{$h \gets Q_{H-1}$ {\bf to} $Q_H-1$}
            \State $r_h \gets \left\lceil \dfrac{(Q_H-h)P_{H-1}+(h-Q_{H-1})P_{H}}{Q_H-Q_{H-1}} - 1 \right\rceil$
            \State $w(r_h,h) \gets \dfrac{b_H r_h+ a_Hh}{b_HP_H+a_HQ_H}$
            \State $s(h) \gets H$
        \EndFor
    \EndFor
    \For{$h \gets Q^*_k+1$ {\bf to} $Q_{k+1}$}
        \State $r_h \gets \left\lceil \dfrac{(h-Q^*_k)P_{k+1}}{Q_{k+1}-Q^*_{k}} - 1 \right\rceil$
        \State $w(r_h,h) \gets \dfrac{b_{k+1} r_h+ a_{k+1}(n-1-h)}{b_{k+1}P_{k+1}+a_{k+1}(n-1-Q_{k+1})}$
        \State $s(h) \gets k+1$
    \EndFor
    \For{$H \gets k+2$ {\bf to} $d$}
        \For{$h \gets Q_{H-1}+1$ {\bf to} $Q_H$}
            \State $r_h \gets \left\lceil \dfrac{(Q_H-h)P_{H-1}+(h-Q_{H-1})P_{H}}{Q_H-Q_{H-1}} - 1 \right\rceil$
            \State $w(r_h,h) \gets \dfrac{b_H r_h+ a_H(n-1-h)}{b_HP_H+a_H(n-1-Q_H)}$
            \State $s(h) \gets H$
        \EndFor
    \EndFor
    \If{$(0,1) \in D$} $s(Q^*_k)\gets0$ \EndIf
\end{algorithmic}
{\bf Step 6:} Ordering the points.
\begin{algorithmic}
    \State Order the points $(r_h,h)$ for $h \gets 0,1, \dots, Q_k-1, Q^*_k+1, Q^*_k+2, \dots, n-1$ after increasing values of $w(r_h,h)$ and call these points in this order $(p_0,q_0), (p_1,q_1) \dots, (p_{N-1},q_{N-1})$.
    \If{$(0,1) \in D$} $(p_{N-1},q_{N-1}) \gets (0,Q^*_k)$ \EndIf
\end{algorithmic}
{\bf Step 7:} Assignment of $f$-values.
\begin{algorithmic}
    \For{$t \gets 1- \max(P_0,P_d)$ {\bf to} $m-M-1$}
        \For{$h \gets 0$ {\bf to} $N-1$}
            \If{$0 \leq p_h+t < m$ \textbf{and} $\delta(p_h+t,q_h) = 0$}
                \State $g(p_h+t,q_h) \gets \lambda(\mathbf{d}_{s(q_h)},p_h+t,q_h)$
                \ForAll{$\mathbf{d} \in D$}
                    \State $\lambda(\mathbf{d},p_h+t,q_h) \gets \lambda(\mathbf{d},p_h+t,q_h)-g(p_h+t,q_h)$
                \EndFor
                \State $\delta(p_h+t,q_h) \gets 1$
            \EndIf
        \EndFor
    \EndFor
\end{algorithmic}
{\bf Step 8:} Start nonvalid case, initial values, cf.~Step 1.
\begin{algorithmic}
    \State $\left((a_1,b_1), \dots, (a_k,b_k)\right) \gets \left((a_k,b_k), \dots, (a_1,b_1)\right)$
    \State $\left((a_{k+1},b_{k+1}), \dots, (a_d,b_d)\right) \gets \left((a_d,b_d), \dots, (a_{k+1},b_{k+1})\right)$
    \If{$(1,0) \in D$} $\mathbf{d}_0 \gets (1,0)$ \EndIf
\end{algorithmic}
{\bf Step 9:} Border points, cf.~Step 2.
\begin{algorithmic}
    \If {$M>m$} $M \gets m$ \EndIf
    \State $\left((P_0,Q_0), \dots, (P_k,Q_k)\right) \gets \left((P_k,Q_k), \dots, (P_0,Q_0)\right)$
    \State $\left((P^*_k, Q^*_k),(P_{k+1},Q_{k+1}), \dots, (P_{d},Q_{d})\right) \gets \left((M-P_{d}-1,n-Q_{d}-1), \dots,\right.$
    \State \hspace*{\fill} $\left.(M-P_{k+1}-1,n-Q_{k+1}-1), (M-P^*_k-1, n-Q^*_k-1)\right)$
\end{algorithmic}
{\bf Step 10:} Choosing starting points $r_h$, weights $w(r_h,h)$ and directions $\mathbf{d}_{s(h)}$, cf.~Step 5.
\begin{algorithmic}
    \For{$H \gets 1$ {\bf to} $k$}
        \For{$h \gets P_{H-1}$ {\bf to} $P_H-1$}
            \State $r_h \gets \left\lceil \dfrac{(P_H-h)Q_{H-1}+(h-P_{H-1})Q_{H}}{P_H-P_{H-1}} - 1 \right\rceil$
            \State $w(h,r_h) \gets \dfrac{a_H r_h+ b_Hh}{a_HQ_H+b_HP_H}$
            \State $s(h) \gets H$
        \EndFor
    \EndFor
    \For{$h \gets P^*_k+1$ {\bf to} $P_{k+1}$}
        \State $r_h \gets \left\lceil \dfrac{(h-P^*_k)Q_{k+1}}{P_{k+1}-P^*_k} - 1 \right\rceil$
        \State $w(h,r_h) \gets \dfrac{M-1-h}{M-1-P^*_k}$
        \State $s(h) \gets k+1$
    \EndFor
    \For{$H \gets k+2$ {\bf to} $d$}
        \For{$h \gets P_{H-1}+1$ {\bf to} $P_H$}
            \State $r_h \gets \left\lceil \dfrac{(P_H-h)Q_{H-1}+(h-P_{H-1})Q_{H}}{P_H-P_{H-1}} - 1 \right\rceil$
            \State $w(h,r_h) \gets \dfrac{a_{H} r_h+ b_{H} (M-1-h)} {a_{H}Q_{H-1} +b_{H}(M-1-P_{H-1})}$
            \State $s(h) \gets H$
        \EndFor
    \EndFor
    \If{$(1,0) \in D$} $s(P^*_k)\gets0$ \EndIf
\end{algorithmic}
{\bf Step 11:} Ordering the points, cf.~Step 6.
\begin{algorithmic}
    \State Order the points $(r_h,h)$ for $h \gets 0,1, \dots, P_k-1, P^*_k+1, P^*_k+2, \dots, M-1$ after increasing values of $w(r_h,h)$ and call these points in this order $(p_0,q_0), (p_1,q_1), \dots, (p_{M-1},q_{M-1})$.
    \If{$(1,0) \in D$} $(p_{M-1},q_{M-1}) \gets (P^*_k,0)$ \EndIf
\end{algorithmic}
{\bf Step 12:} Assignment of $f$-values, cf.~Step 7.
\begin{algorithmic}
    \State $u \gets m-M$
    \For{$t \gets 1- \max(Q_0,Q_d)$ {\bf to} $n-1$}
        \For{$h \gets 0$ {\bf to} $M-1$}
            \If{$0 \leq p_h+u < m$ \textbf{and} $0 \leq q_h+t < n$ \textbf{and} $\delta(p_h+u,q_h+t) = 0$}
                \State $g(p_h+u,q_h+t) \gets \lambda(\mathbf{d}_{s(p_h)},p_h+u,q_h+t)$
                \ForAll{$\mathbf{d} \in D$}
                    \State $\lambda(\mathbf{d},p_h+u,q_h+t) \gets \lambda(\mathbf{d},p_h+u,q_h+t)-g(p_h+u,q_h+t)$
                \EndFor
                \State $\delta(p_h+u,q_h+t) \gets 1$
            \EndIf
        \EndFor
    \EndFor
    \Return{g}
\end{algorithmic}

\begin{rem}{\rm
We are assuming that the line sums are exact. So, the fact that the output is consistent with the data can be easily checked by observing whether all the line sums are equal to zero (Steps 7 and 12 update the line sums by subtracting the value of each point).

}
\end{rem}

\begin{ex}\label{ex_steps}{\rm
Let be given $m=21$, $n=16$, $D=\left\{(0,1),(1,0),(1,1),(-1,1)\right.$, $\left.(-3,-1),(-1,-3),(5,-1),(7,5)\right\}$ and the line sums in the directions of $D$ of some function $f:A \to \R$ (the function itself is irrelevant for the example). The effect of the steps is the following.
\vskip.1cm
\noindent {\bf Step 1.} Initial values: $k=2,d=6$, $\mathbf{d}_0=(0,1)$, $(a_1,b_1)=(5,1)$, $(a_2,b_2)=(1,1)$, $(a_3,b_3)=(1,3)$, $(a_4,b_4)=(1,1)$, $(a_5,b_5)=(7,5)$, $(a_6,b_6)=(3,1)$. $M=19, N=13$.
\vskip.1cm
\noindent {\bf Step 2.} Border points: $(P_0,Q_0)=(6,0)$, $(P_1,Q_1)=(1,1)$, $(P_2,Q_2)=(0,2)$, $(P^*_2,Q^*_2)=(0,5)$, $(P_3,Q_3)=(1,8)$, $(P_4,Q_4)=(2,9)$, $(P_5,Q_5)=(9,14)$, $(P_6,Q_6)=(12,15)$.
\vskip.1cm
\noindent {\bf Step 3} Switching functions: $\delta(p,q)=1, g(p,q)=0$ for $p=0,1$ and $q=2,3,4$; $\delta(p,q)=0$ for all other points $(p,q)$ of $A$.
\vskip.1cm
\noindent {\bf Step 4} False: We are in a valid case, since $m > M$.

\begin{figure}[htbp]
\centering
\includegraphics[scale=0.5]{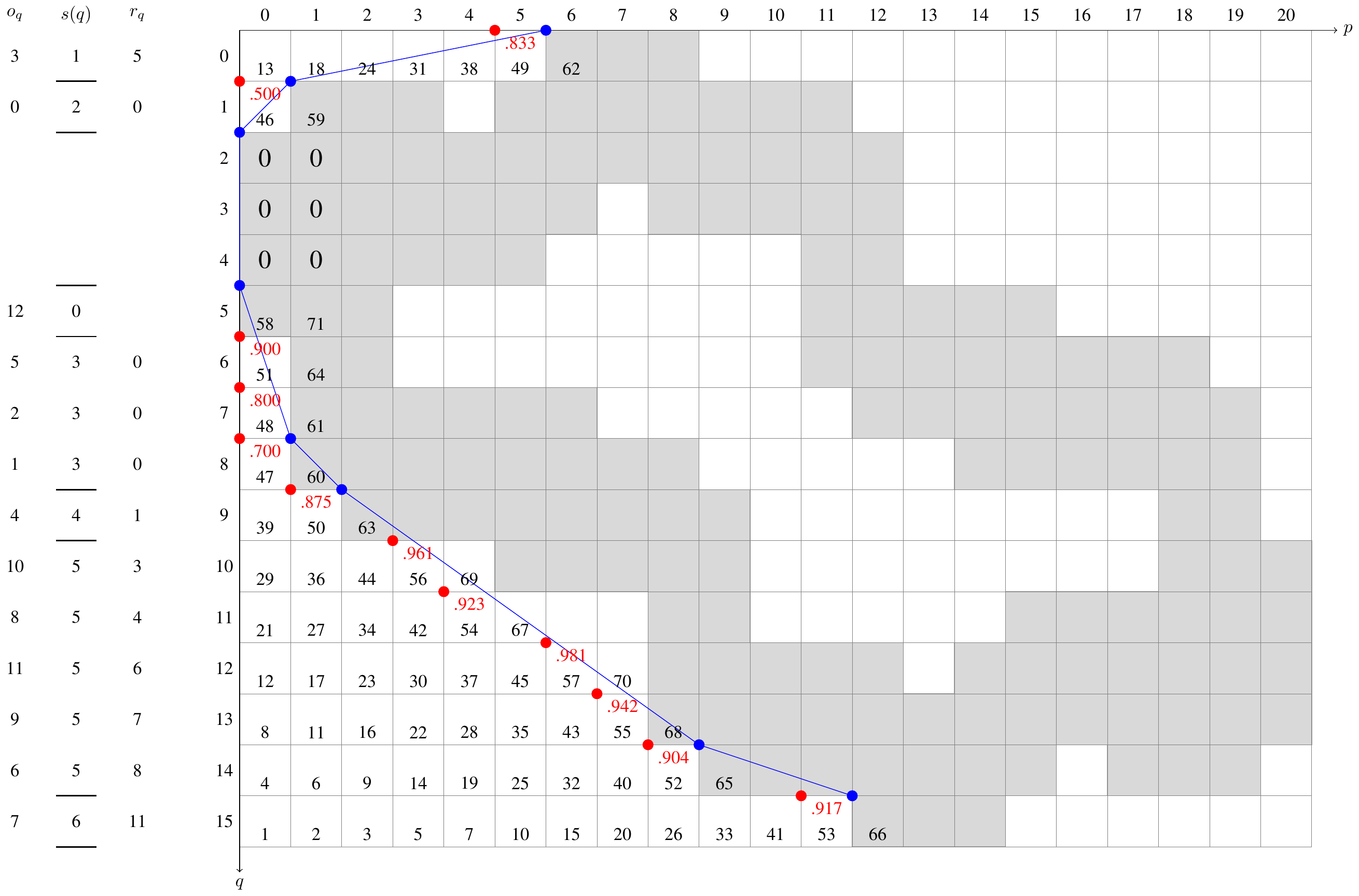}
\caption{\small An illustration of Steps 1-7 of Example \ref{ex_steps}. The light grey pixels indicate the union of the switching domains. The white pixels indicate the pixels of which the $f$-values are unique, and therefore equal to the computed $g$-value. There are six primitive switching functions and their lexicographic smallest elements have a 0. In Step 3 their $g$-values are fixed as $0$, but this may be replaced by any other values. Step 1 guarantees the concavity of the upper left corner region $V_{UL}$ and the lower left corner region $V_{LL}$, left of the broken line. The border points for $V_{UL}$ and $V_{LL}$ are indicated by the dots along the broken line. They are found in Step 2. Step 4 provides a shortcut in case of a nonvalid case; if $m \leq M$, then the coordinates can be switched. For each row the grid point just left of the border line is computed in Step 5. The weights of these points are indicated in the upper numbers inside the pixels. The (highlighted) points immediately left of the broken line are ordered after size as indicated in the column $o_q$ (see Step 6). The function $s$ indicates the directions which are used for the grid points in that row. Finally, in Step 7, the $g$-values are computed for the first $m-M$ columns and some more pixels. The order in which they are calculated is given in black.}
\label{fig:algvalid}
\end{figure}

\vskip.1cm
\noindent {\bf Step 5.} Choice of starting points, weights and directions: $r_0=5$, $r_1=0$, $r_6=r_7=r_8=0$, $r_9=1$, $r_{10}=3$, $r_{11}=4$, $r_{12}=6$, $r_{13}=7$, $r_{14}=8$, $r_{15}=11$; $w(5,0)=.833$, $w(0,1)=.500$, $w(0,6)=.900$, $w(0,7)=.800$, $w(0,8)=.700$, $w(1,9)=.875$, $w(3,10)=.962$, $w(4,11)=.923$, $w(6,12)=.981$, $w(7,13)=.942$, $w(8,14)=.904$, $w(11,15)=.917$; $s(0)=1$, $s(1)=2$, $s(6)=s(7)=s(8)=3$, $s(9)=4$, $s(10)=s(11)=s(12)=s(13)=s(14)=5$, $s(15)=6$, $s(5)=0$ (See Figure \ref{fig:algvalid}).
\vskip.1cm
\noindent {\bf Step 6.} Ordering of the points: $(p_0,q_0)=(0,1)$, $(p_1,q_1)=(0,8)$, $(p_2,q_2)=(0,7)$, $(p_3,q_3)=(5,0)$, $(p_4,q_4)=(1,9)$, $(p_5,q_5)=(0,6)$, $(p_6,q_6)=(8,14)$, $(p_7,q_7)= (11,15)$, $(p_8,q_8)=(4,11)$, $(p_9,q_9)=(7,13)$, $(p_{10},q_{10})=(3,10)$, $(p_{11},q_{11})=(6,12)$, $(p_{12},q_{12})=(0,5)$.
\vskip.1cm
\noindent {\bf Step 7.} Assignment. See Figure \ref{fig:algvalid} for the order in which the $f$-values are computed, indicated by the black numbers. After this step the $f$-values of the first $m-M=2$ columns are known and $M=19$ columns are left.
\vskip.1cm
\noindent {\bf Step 8.} We are now in a nonvalid case and apply a switch of coordinate axes. The new initial values are: $\mathbf{d}_0=(1,0)$, $(a_1,b_1)=(1,1)$, $(a_2,b_2)=(5,1)$, $(a_3,b_3)=(3,1)$, $(a_4,b_4)=(7,5)$, $(a_5,b_5)=(1,1)$, $(a_6,b_6)=(1,3)$.
\vskip.1cm
\noindent {\bf Step 9.} Border points: $(P_0,Q_0)=(0,2)$, $(P_1,Q_1)=(1,1)$, $(P_2,Q_2)=(6,0)$, $(P^*_2,Q^*_2)=(6,0)$, $(P_3,Q_3)=(9,1)$, $(P_4,Q_4)=(16,6)$, $(P_5,Q_5)=(17,7)$, $(P_6,Q_6)=(18,10)$.
\vskip.1cm
\noindent {\bf Step 10.} Choice of starting points, weights and directions: $r_0=1, r_1=r_2= \dots = r_5=0$, $r_7=r_{8}=r_{9}=0$, $r_{10}=1$, $r_{11}=2$, $r_{12}=r_{13}=3$, $r_{14}=4$, $r_{15}=r_{16}=5$, $r_{17}=6$, $r_{18}=9$; $w(0,1)=.500$, $w(1,0)=.167$, $w(2,0)=.333$, $w(3,0)=.500$, $w(4,0)=.667$, $w(5,0)=.833$, $w(7,0)=.917$, $w(8,0)=.833$, $w(9,0)=.750$, $w(10,1)=.904$, $w(11,2)=.942$, $w(12,3)=.981$, $w(13,3)=.885$, $w(14,4)=.923$, $w(15,6)=.962$, $w(16,5)=.865$, $w(17,6)=.875$, $w(18,9)=.900$; $s(0)=1, s(1)=s(2)=\dots =s(5)=2$, $s(7)=s(8)=s(9)=3$, $s(10)=s(11)=\dots =s(16)=4$, $s(17)=5$, $s(18)=6$, $s(6)=0$.
\vskip.1cm
\noindent {\bf Step 11} Ordering of the points: $(p_0,q_0)=(1,0)$, $(p_1,q_1)=(2,0)$, $(p_2,q_2)=(0,1)$, $(p_3,q_3)=(3,0)$, $(p_4,q_4)=(4,0)$, $(p_5,q_5)=(9,0)$, $(p_6,q_6)=(5,0)$, $(p_7,q_7)= (8,0)$, $(p_8,q_8)=(16,5)$, $(p_9,q_9)=(17,6)$, $(p_{10},q_{10})=(13,3)$, $(p_{11},q_{11})=(18,9)$, $(p_{12},q_{12})=(10,1)$, $(p_{13},q_{13})=(7,0)$, $(p_{14},q_{14})=(14,4)$, $(p_{15},q_{15})=(11,2)$, $(p_{16},q_{16})=(15,5)$, $(p_{17},q_{17})=(12,3)$, $(p_{18},q_{18})=(6,0)$.
\vskip.1cm
\noindent {\bf Step 12.} Assignment. See Figure \ref{fig:algnonvalid} for the order in which the $g$-values are computed, indicated by the black numbers. This has to be continued in the obvious way. When this step has been completed all the $g$-values are known.}
\end{ex}

\section{Complexity}\label{sec:complexity}
We state the complexity of each step of the algorithm, where we count every addition, subtraction, multiplication, division and determining of the larger of two explicit quantities as one operation.
\begin{itemize}
\item[-] Step 1: $\mathcal{O}(d \log d)$;
\item[-] Step 2: $\mathcal{O}(d)$;
\item[-] Step 3: $\mathcal{O}(mn)$;
\item[-] Step 4: $\mathcal{O}(1)$;
\item[-] Steps 5 and 10: $\mathcal{O}(n)$ and $\mathcal{O}(m)$, respectively;
\item[-] Steps 6 and 11: $\mathcal{O}(n \log n)$ and $\mathcal{O}(m \log m)$, respectively;
\item[-] Steps 7+12: $\mathcal{O}(dmn)$.
\end{itemize}

Without loss of generality we may assume $m \leq n$. It follows that the complexity of the algorithm is $\mathcal{O}(dmn)$ unless
$mn = \mathcal{O}(\log d)$ or $dm = \mathcal{O}(\log n)$ or $dn = \mathcal{O}(\log m)$. The latter case can not occur since $m \leq n$.
If $m = \mathcal{O}(\log d)$, then we recall that $a_h \geq 1$ with at most one exception. So we can delete $d-m-1$ directions and still have a nonvalid case. After deletion we have, denoting by $d$ the new number of directions, $m = d-1$ and complexity $\mathcal{O}(dmn)$.

It remains to consider $dm = \mathcal{O}(\log n)$. In this case only the complexity of Step 6 has to be adjusted. This can be achieved by remembering in Step 5 for every $H$ the location of the point $(r_h,h)$ with $s(h)=H$ and minimal weight for direction $(a_H,b_H)$. This has complexity $\mathcal{O}(dn)$.
Now Step 6 proceeds as follows. For each direction $(a_H,b_H)$ let $(r_h,h)$ be the point with $s(h)=H$ and minimal weight. For $H\leq k$ we start with the vectors
$$(r_h,\ h,\ b_H(P_H-1)+a_HQ_H,\ b_HP_H+a_HQ_H,\ R_H,\ S_H,\ T_H,\ U_H),$$
where $(R_H,S_H)=(Q_{H-1},Q_H-1)$ and $(T_H,U_H)$ is the unique pair satisfying $0 \leq T_H<a_H$ and $b_HT_H+a_HU_H = 1$. For $k<H\leq d$ we start with the vectors
$$(r_h,\, h,\ b_H(P_H-1)+a_H(n-1-Q_H),\ b_HP_H+a_H(n-1-Q_H),\ R_H, S_H,T_H, U_H),$$
where $(R_H,S_H)=(Q_{H-1}+1,Q_H)$ and $(T_H,U_H)$ is the unique pair satisfying $0 \leq T_H<a_H$ and $b_HT_H-a_HU_H = 1$. Observe that in each case the quotient of the third and fourth entry is the weight.
We order such vectors on the top line according to increasing weight. At each step we increase by one the third entry of the first (leftmost) vector. If the third entry now is still smaller than the fourth entry, we replace the first two entries $(r_h,h)$ by $(r_h+T_H,h+U_H)$ or $(r_h+T_H-b_H,h+U_H-a_H)$ such that the second entry is in $[R_H,S_H]$. If the third entry becomes equal to the fourth one, we neglect the vector in the sequel. In any case we order the remaining vectors on the line again after increasing weight. This procedure runs until there is no vector left. At every step the two leftmost entries of the leftmost vector give the next value $(p_h,q_h)$.
The computation of the vectors $(T_H,U_H)$ has complexity $\mathcal{O}(d \log n)$, the computation of each row $\mathcal{O}(d \log d)$ and there are $n$ rows. Therefore the total complexity is $\mathcal{O}(nd \log d)$. This is $\mathcal{O}(dmn)$, unless $m = \mathcal{O}(\log d)$. We have already remarked that in this case  we can delete $d-m-1$ directions, still have a nonvalid case, and have complexity $\mathcal{O}(dmn)$.

\begin{ex}\label{ex_stepsctd} {\rm [Continuation of Example \ref{ex_steps}].
The described procedure yields as the first row the vectors
\begin{displaymath}
\begin{array}{lll}
(0,1,1,2,1,1,0,1), & (0,8,7,10,6,8,0,-1), & (5,0,5,6,0,0,1,0),\\
(1,9,7,8,9,9,1,0), & (8,14,47,52,10,14,3,2), & (11,15,11,12,15,15,1,0).
\end{array}
\end{displaymath}
Since the last four entries do not change, we do not mention them in the table. Then the table becomes as follows (each row represents one step in the procedure).
\begin{center}
\footnotesize{
\setlength{\tabcolsep}{3pt}
\begin{tabular}{ccccccc}
(0,1,1,2) & (0,8,7,10) & (5,0,5,6) & (1,9,7,8) & (8,14,47,52) &(11,15,11,12) \\
(0,8,7,10) & (5,0,5,6) & (1,9,7,8) & (8,14,47,52) & (11,15,11,12)& \\
(0,7,8,10) & (5,0,5,6) & (1,9,7,8) & (8,14,47,52) &(11,15,11,12)& \\
(5,0,5,6) & (1,9,7,8) & (0,6,9,10) & (8,14,47,52)& (11,15,11,12)& \\
(1,9,7,8) & (0,6,9,10) & (8,14,47,52) &(11,15,11,12)&& \\
(0,6,9,10) & (8,14,47,52)& (11,15,11,12)& &&\\
(8,14,47,52)& (11,15,11,12)&&&&\\
(11,15,11,12)& (4,11,48,52)& &&&\\
(4,11,48,52)&&&&&\\
(7,13,49,52)&&&&&\\
(3,10,50,52)&&&&&\\
(6,12,51,52)&&&&&\\
\end{tabular}}
\end{center}

\noindent The sequence $(p_0,q_0)=(0,1)$, $(p_1,q_1)=(0,8)$, $(p_2,q_2)=(0,7)$, \dots, $(p_{11},q_{11})=(6,12)$ can be read from the leftmost two entries. At the end $(p_{12},q_{12})=(0,5)$ has to be added.
}
\end{ex}



\section{Conclusions}\label{sec:conclusion}
In this paper we have addressed the tomographic reconstruction problem for functions with values in a unique factorization domain or field, such as integers and reals. A key argument is that one may ask for the point values even when many solutions are admissible, since the values of points outside the switching domains are common to all functions satisfying the problem.

Starting from the characterization of the switching functions in \cite{hatij} and the results in \cite{pati}, we have shown that all points with uniquely determined value, namely, not belonging to switching domains, can be recovered once we give an arbitrary value to $(m-M)(n-N)$ points, where $(m-M)(n-N)$ is the number of linearly independent switching functions. We have provided an algorithm which computes the point values systematically and runs in time linear in $dmn$, where $d$ is the number of directions. By the result in \cite{hatij} our algorithm provides the complete set of solutions with values in the unique factorization domain or in the field.

The proposed approach works when line sums are supposed to be exact and therefore not all projections are necessary to recover a solution. It underlines the structural difference between unique factorization domains and other kinds of sets, such as $\{0,1\}$ (leading to binary images), since in the latter case the reconstruction problem has proven to be NP-hard \cite{ggp99,ggp00}.

Two questions arise by this paper. Firstly, does there exist a similar algorithm for higher dimensions? Secondly, given a system of inconsistent line sums, is there a fast way to find a best approximation of consistent line sums so that the algorithm of this paper can be applied to construct the most likely set of solutions (over $\Z$ or $\R$)?

\section*{References}
\bibliographystyle{plain}
\bibliography{cst_references}

\end{document}